\documentclass{amsart}
\usepackage{amsfonts}
\usepackage{amsmath}
\usepackage{enumerate}
\usepackage{amsthm}
\usepackage{hyperref}
\usepackage{cite}

\begin{document}

\newtheorem{definition}{Definition}[section]
\newtheorem{theorem}[definition]{Theorem}
\newtheorem{proposition}[definition]{Proposition}
\newtheorem{remark}[definition]{Remark}
\newtheorem{lemma}[definition]{Lemma}
\newtheorem{corollary}[definition]{Corollary}
\newtheorem{example}[definition]{Example}

\numberwithin{equation}{section}

\title[Tamed exhaustion functions and Schwarz Lemma]{Tamed exhaustion functions and Schwarz type lemmas for almost Hermitian manifolds}

\author{Weike Yu}

\date{}

\begin{abstract}
In this paper, we study a special exhaustion function on almost Hermitian manifolds and establish the existence result by using the Hessian comparison theorem. From the viewpoint of the exhaustion function, we establish related Schwarz type lemmas for almost holomorphic maps between two almost Hermitian manifolds. As corollaries, we deduce Liouville type theorems for almost holomorphic maps.
\end{abstract}
\keywords{Almost Hermitian manifold; Canonical connection; Hessian comparison theorem; Tamed exhaustion function; Schwarz lemma.}
\subjclass[2010]{32Q60}
\maketitle
\section{Introduction}
The classical Schwarz lemma which was reformulated by Pick states that a holomorphic map between unit discs $D$ in $\mathbb{C}$ decreases the Poincar\'e metric. Later, Ahlfors \cite{[Ah]} generalized this lemma to holomorphic maps from the disc $D$ into a negative curved Riemannian surface. Schwarz lemma is extremely useful in complex analysis and differential geometry, and has been extended to several cases, such as holomorphic maps between higher dimensional complex manifolds (cf. \cite{[Ch], [Lu], [Yau], [Ro], [Ni]}, etc.), conformal immersions and harmonic maps between Riemannian manifolds (cf. \cite{[RRV],[Yau-1],[GI], [GIP], [Sh]}, etc.), generalized holomorphic maps between CR manifolds and Hermitian maniflolds (cf. \cite{[DRY], [CDRY], [RT]}, etc.), almost holomorphic maps between almost Hermitian manifolds (cf. \cite{[To]}) and so on.

The precent paper is devoted to establishing Schwarz type lemmas for almost holomorphic maps between almost Hermitian manifolds from the viewpoint of the tamed exhaustion function (see Definition \ref{def3.2}), which was introduced by Royden in \cite{[Ro]}. Roughly speaking, an almost Hermitian manifold $(M,J,g)$ is an almost complex manifold $(M,J)$ equipped with a Riemannian metric $g$ compatible with the almost complex structure $J$. On an almost Hermitian manifold $(M,J,g)$, there is a preferred connection preserving the metric $g$ and the almost complex structure $J$, and coinciding with the Chern connection when $J$ is integrable. Such a connection is called the canonical connection, which was first introduced by Ehresman and Libermann (cf. \cite{[EL]}). Although it has nontrivial torsion, it is more suitable for analytic problem on almost Hermitian manifolds than the Levi-Civita connection. So we will always use the canonical connection in present paper unless otherwise specified. Firstly, we will establish the existence of the tamed exhaustion function on almost Hermitian manifolds by using Hessian comparison theorem.
\begin{theorem}\label{thm1.1}
Let $(M^m,J, g)$ be a complete almost Hermitian manifold with holomorphic bisectional curvature 
\begin{align}
R(X,\overline{X},Y,\overline{Y})\geq -B(1+r(x))^2,\label{1.1}
\end{align}
the torsion 
\begin{align}
\|\tau(X,Y)\|\leq A_1(1+r(x)),\label{1.2}
\end{align}
and the $(2,0)$ part of the curvature tensor 
\begin{align}
|R(\overline{X},Y,Y,X)|\leq A_2(1+r(x))^2,\label{1.3}
\end{align}
where $r(x)$ is the Riemannian distance of $M$ between a fixed point $x_0$ and $x$ in $M$, $X,Y\in T^{1,0}_{x}M$ with $\|X\|=\|Y\|=1$, and $B, A_1, A_2$ are positive constants. Then there is a tamed exhaustion function on $M$ in the sense of Definition \ref{def3.2}.
\end{theorem}

Next, we will establish a Schwarz type lemma for almost holomorphic maps from an almost Hermitian manifold admitting a tamed exhaustion function into a general almost Hermitian manifold, which is a generalization of the result in \cite[Section 3]{[Ro]}.
\begin{theorem}\label{thm1.2}
Let $(M^m,J, g)$ be an almost Hermitian manifold with holomorphic sectional curvature bounded below by $-k_1\ (k_1\geq 0)$. Let $(\tilde{M}^n,\tilde{J},\tilde{g})$ be an almost Hermitian manifold with holomorphic sectional curvature bounded above by $-k_2\ (k_2>0)$. Assume that $(M^m,J, g)$ admits a tamed exhaustion function. Then for any almost holomorphic map $f: M\rightarrow \tilde{M}$, we have
\begin{align}
f^*\tilde{g}\leq \frac{k_1}{k_2}g.
\end{align}
In particular, if $k_1=0$, then every almost holomorphic map is constant.
\end{theorem}
Combining the Theorem \ref{thm1.1} and Theorem \ref{thm1.2}, we obtain another version of Schwarz type lemma where assumptions on domain and target manifolds are all about curvatures and torsions. 
\begin{corollary}
Let $(M^m,J, g)$ be a complete almost Hermitian manifold with holomorphic sectional curvature bounded below by $-k_1\ (k_1\geq 0)$. Let $(\tilde{M}^n,\tilde{J},\tilde{g})$ be an almost Hermitian manifold with holomorphic sectional curvature bounded above by $-k_2\ (k_2>0)$. Assume that $(M^m,J, g)$ satisfies \eqref{1.1}, \eqref{1.2} and \eqref{1.3}. Then for any almost holomorphic map $f: M\rightarrow \tilde{M}$, we have
\begin{align}
f^*\tilde{g}\leq \frac{k_1}{k_2}g.
\end{align}
In particular, if $k_1=0$, then every almost holomorphic map is constant.
\end{corollary}
Note that the above corollary is an extension of \cite[Theorem 1]{[YC]} in Hermitian case and \cite{[CCL]}, \cite[Theorem 2]{[Ro]}, \cite[Theorem 1.4]{[Ni]} in K\"ahler case. 

We point out that when $M$ is compact, the assumptions \eqref{1.1}, \eqref{1.2}, \eqref{1.3} are automatically satisfied. Hence, we obtain the following
\begin{corollary}\label{cor1.5}
Let $(M^m,J, g)$ be a compact almost Hermitian manifold with holomorphic sectional curvature bounded below by $-k_1\ (k_1\geq 0)$. Let $(\tilde{M}^n,\tilde{J},\tilde{g})$ be an almost Hermitian manifold with holomorphic sectional curvature bounded above by $-k_2\ (k_2>0)$. Then for any almost holomorphic map $f: M\rightarrow \tilde{M}$, we have
\begin{align}
f^*\tilde{g}\leq \frac{k_1}{k_2}g.
\end{align}
\end{corollary}
In this case, we also have a Liouville type theorem as follows.
\begin{corollary}\label{cor1.6}
Let $(M^m,J, g)$ be a compact almost Hermitian manifold with nonnegative (resp. positive) holomorphic sectional curvature. Let $(\tilde{M}^n,\tilde{J},\tilde{g})$ be an almost Hermitian manifold with negative (resp. nonpositive) holomorphic sectional curvature. Then any almost holomorphic map $f: M\rightarrow \tilde{M}$ is constant.
\end{corollary}
We remark that the above corollary has been proved in \cite{[Ma]} by a different method, which is an extension of \cite[Theorem 1.2]{[Ya]} to the almost Hermitian case.

\section{Preliminaries}\label{section2}
In this section, we will introduce some notions and notations of almost Hermitian geometry (cf. \cite{[To], [TWY]}).

Let $(M^m,J)$ be an almost complex manifold with $\dim_{\mathbb{C}} M=m$. An almost Hermitian metric $g$ on $(M,J)$ is a Riemannian metric with 
\begin{align}
g(JX,JY)=g(X,Y)
\end{align}
for any tangent vectors $X,Y\in TM$. The triple $(M,J,g)$ is referred as to an almost Hermitian manifold. Let $TM^{_\mathbb{C}}=TM\otimes_{\mathbb{R}} \mathbb{C}$ denote the complexified tangent space of $M$, and we extend the almost complex structure $J$ and the almost Hermitian metric $g$ from $TM$ to $TM^{\mathbb{C}}$ by $\mathbb{C}-$linearity. Let $T^{1,0}M$ (resp. $T^{0,1}M$) denote eigenspace of $J$ corresponding to the eigenvalue $\sqrt{-1}$ (resp. $-\sqrt{-1}$), then one has the following decomposition:
\begin{align}
TM^{\mathbb{C}}=T^{1,0}M\oplus T^{0,1}M,
\end{align}
where $T^{1,0}M=\{X-\sqrt{-1}JX: X\in TM\}$ and $T^{0,1}M=\overline{T^{1,0}M}$. We note that by extending the almost complex structure $J$ to forms, every $m$-form can be decomposed into $(p,q)$-forms for each $p,q\geq 0$ with $p+q=m$.

Let $\nabla$ be an affine connection on $TM$, which we extend $\mathbb{C}-$linearly to $TM^{\mathbb{C}}$.
This connection $\nabla$ is called an almost Hermitian connection if it satisfies
\begin{align}
\nabla g=0,\ \nabla J=0.
\end{align}
Analogue to Riemannian geometry, the torsion of the connection $\nabla$ is defined by
\begin{align}
\tau(X,Y)=\nabla_XY-\nabla_YX-[X,Y].
\end{align}
Note that An almost Hermitian connection is uniquely determined by the $(1,1)$-part of the torsion $\tau$ (cf. \cite{[Ga]}). In particular, there is a unique almost Hermitian connection $\nabla$ on $(M,J,g)$ whose torsion $\tau$ has vanishing $(1,1)$-part everywhere (see \cite{[Ko-1]}). Such a connection is usually called the canonical connection of the almost Hermitian manifold.

Suppose $(M^m,J,g)$ is an almost Hermitian manifold with canonical connection $\nabla$. Let $\{e_i\}_{i=1}^m$ be a local unitary frame field of $T^{1,0}M$, and $\{\theta^i\}_{i=1}^m$ be its coframe field. Then the almost Hermitian metric $g$ can be expressed as
\begin{align}
g=\sum_{i,k}g_{j\bar{k}}\theta^j\theta^{\bar{k}}=\sum_j\theta^j\theta^{\bar{j}},
\end{align}
and the torsion $\tau$ may be written as 
\begin{align}
\tau=\sum_i(\tau^ie_i+\tau^{\bar{i}}e_{\bar{i}}),
\end{align}
where
\begin{align}
\tau^i=\frac{1}{2}\sum_{j,k}\left(\tau^i_{jk}\theta^j\wedge\theta^k+\tau^i_{\bar{j}\bar{k}}\theta^{\bar{j}}\wedge\theta^{\bar{k}}\right)
\end{align}
with $\tau^i_{jk}=-\tau^i_{kj}$ and $\tau^i_{\bar{j}\bar{k}}=-\tau^i_{\bar{k}\bar{j}}$. Note that the $(0,2)$-part of $\tau^i$ is independent of the choice of metrics, and indeed can be regarded as the Nijenhuis tensor of the almost complex structure $J$ (see, e.g., \cite{[TWY]}).

According to \cite{[To]}, we have the structure equations for the canonical connection $\nabla$ as follows:
\begin{align}
d\theta^i=-\sum_j\theta^i_j\wedge\theta^j+\tau^i,
\end{align}
\begin{align}
d\theta^i_j=-\sum_k\theta^i_k\wedge\theta^k_j+\Omega^i_j,
\end{align}
with
\begin{align}
\theta^i_j+\theta^{\bar{j}}_{\bar{i}}=0,
\end{align}
and 
\begin{align}
\Omega^i_j=\sum_{k,l}\left(\frac{1}{2}R^i_{jkl}\theta^k\wedge\theta^l+R^i_{jk\bar{l}}\theta^k\wedge\theta^{\bar{l}}+\frac{1}{2}R^i_{j\bar{k}\bar{l}}\theta^{\bar{k}}\wedge\theta^{\bar{l}}\right),
\end{align}
where $\theta^i_j$ is connection $1$-form of $\nabla$, and $R^i_{jkl}, R^i_{jk\bar{l}}, R^i_{j\bar{k}\bar{l}}$ are components of the curvature tensor of $\nabla$, that is
\begin{align}
R(X,Y)Z=\nabla_X\nabla_YZ-\nabla_Y\nabla_XZ-\nabla_{[X,Y]}Z,
\end{align}
\begin{align}
R(X,Y,Z,W)=g(R(Z,W)X,Y).
\end{align}
Set 
\begin{align}
R_{j\bar{i}kl}=R^i_{jkl},\ R_{j\bar{i}k\bar{l}}=R^i_{jk\bar{l}},\ R_{j\bar{i}\bar{k}\bar{l}}=R^i_{j\bar{k}\bar{l}},
\end{align}
then one have the following properties (cf. \cite{[TWY]}):

\begin{equation}
R_{j\bar{i}kl}=\overline{R_{i\bar{j}\bar{l}\bar{k}}},\ R_{j\bar{i}k\bar{l}}=\overline{R_{i\bar{j}l\bar{k}}},
\end{equation}
\begin{equation}
R_{j\bar{i}kl}=-R_{\bar{i}jkl}=-R_{j\bar{i}lk},\ R_{j\bar{i}k\bar{l}}=-R_{\bar{i}jk\bar{l}}=-R_{j\bar{i}\bar{l}k},\\
\end{equation}
\begin{align}
R_{i\bar{j}k\bar{l}}-R_{k\bar{j}i\bar{l}}&=\tau^j_{ik;l}-\tau^{\bar{m}}_{ik}\tau^j_{\bar{l}\bar{m}},\\
R_{i\bar{j}kl}&=-\tau^{\bar{i}}_{kl;\bar{j}}+\tau^{\bar{i}}_{\bar{j}\bar{m}}\tau^{\bar{m}}_{kl},\label{2.18.}
\end{align}
where $";"$ means taking covariant derivatives with respect to canonical connection $\nabla$.

Suppose $X=\sum_iX^ie_i$ and $Y=\sum_jY^je_j$ are two $(1,0)$ vectors, then the holomorphic bisectional curvature in direction $X$ and $Y$ is defined as
\begin{align}
\frac{\sum_{i,j,k,l}R_{i\bar{j}k\bar{l}}X^iX^{\bar{j}}Y^kY^{\bar{l}}}{\left(\sum_iX^iX^{\bar{i}}\right)\left(\sum_iY^iY^{\bar{i}}\right)}.
\end{align}
If $X=Y$, the above quantity is called the holomorphic sectional curvature in direction $X$. The first and Second Ricci curvature are respectively defined by
\begin{align}
R'_{k\bar{l}}=\sum_{i}R_{i\bar{i}k\bar{l}},\ R''_{i\bar{j}}=\sum_{k}R_{i\bar{j}k\bar{k}},
\end{align}
thus the scalar curvature is given by
\begin{align}
R=\sum_kR'_{k\bar{k}}=\sum_kR''_{k\bar{k}}.
\end{align}

For a smooth function $u: (M,J,g)\rightarrow \mathbb{R}$, we can define the Hessian of $u$ with respect to the canonical connection $\nabla$ as
\begin{align}
(\nabla du)(X,Y)=Y(Xu)-(\nabla_YX)u\label{2.20}
\end{align}
for any $X,Y\in TM$. Set $u_{k\bar{l}}=(\nabla du)(e_k,e_{\bar{l}})$. According to the Ricci identity (cf. \cite[Lemma 2.1]{[Yu]}), we obtain
\begin{align}
u_{k\bar{l}}=u_{\bar{l}k},
\end{align}
thus the canonical Laplacian of $u$ can be defined by
\begin{align}
\Delta u=trace(\nabla du)=\sum_k\left(u_{k\bar{k}}+u_{\bar{k}k} \right)=2\sum_k\left(u_{k\bar{k}}\right).
\end{align}
On the other hand, let $\nabla^{LC}$ denote the Levi-Civita connection on $(M,J,g)$, then the Hessian of $u$ with respect to $\nabla^{LC}$ can be similarly defined as in \eqref{2.20}. From the relationship between the Levi-Civita connection and the canonical connection (cf. \cite{[Ga]}):
\begin{align}
\left\langle\nabla^{LC}_YX-\nabla_YX,Z\right\rangle=\frac{1}{2}\left(\left\langle\tau(X,Y),Z\right\rangle+\left\langle\tau(Y,Z),X\right\rangle-\left\langle\tau(Z,X),Y\right\rangle \right),
\end{align}
where $\left\langle\cdot,\cdot\right\rangle$ is the inner product induced by the metric $g$, it follows that
\begin{align}
u_{k\bar{l}}-u_{,k\bar{l}}=\left\langle\nabla^{LC}_{e_{\bar{l}}}e_k-\nabla_{e_{\bar{l}}}e_k, \nabla u \right\rangle=\frac{1}{2}\sum_i\left(\tau^l_{ki}u_{\bar{i}}+\tau^{\bar{k}}_{\bar{l}\bar{i}}u_{\bar{i}} \right),
\end{align}
where $u_{,k\bar{l}}=(\nabla^{LC} du)(e_k,e_{\bar{l}})$. Hence, 
\begin{align}
\Delta u-\Delta^{LC}u=\sum_{i,k}\left(\tau^k_{ki}u_{\bar{i}}+\tau^{\bar{k}}_{\bar{k}\bar{i}}u_{\bar{i}} \right),
\end{align}
where $\Delta^{LC}$ denote the Laplacian with respect to $\nabla^{LC}$.

Suppose $(\tilde{M}^n,\tilde{J}, \tilde{g})$ is another almost Hermitian manifold with canonical connection $\tilde{\nabla}$. Let $\{\tilde{e}_\alpha\}_{\alpha=1}^n$ be a unitary frame field of $T^{1,0}\tilde{M}$ on a domain of $\tilde{M}$, and let $\{\tilde{\theta}^\alpha\}_{\alpha=1}^n$ be its coframe field. For simplicity, we will use the same notations as in $M$ for the corresponding geometric data of $\tilde{M}$, such as connection $1$-forms, curvature tensor and torsion, etc., but with $\tilde{}$ on them. 

A smooth map $f: (M,J,g)\rightarrow (\tilde{M},\tilde{J},\tilde{g})$ is said to be almost holomorphic if
\begin{align}
df\circ J=\tilde{J}\circ df.
\end{align}
Since $f$ is almost holomorphic, there are functions $f^\alpha_i$ locally defined on $M$ such that
\begin{align}
f^*\tilde{\theta}^\alpha=\sum_if^\alpha_i\theta^i.\label{2.27}
\end{align}
Taking the exterior derivative of \eqref{2.27} and using the structure equations of $M$ and $\tilde{M}$, we have
\begin{align}
\sum_iDf^\alpha_i\wedge\theta^i+\sum_if^\alpha_i\tau^i-f^*\tilde{\tau}^\alpha=0,\label{2.28}
\end{align}
where
\begin{align}
Df^\alpha_i=df^\alpha_i-\sum_kf^\alpha_k\theta^k_i+\sum_\beta f^\beta_i\hat{\theta}^\alpha_\beta=\sum_k\left(f^\alpha_{ik}\theta^k+f^\alpha_{i\bar{k}}\theta^{\bar{k}}\right).\label{2.29}
\end{align}
Note that for simplicity, we denote $\hat{\theta}^\alpha_\beta=f^*\theta^\alpha_\beta, \hat{\tau}^\alpha_{\beta\gamma}=f^*{\tau}^\alpha_{\beta\gamma},\hat{R}^\alpha_{\beta\gamma\bar{\delta}}=f^*R^\alpha_{\beta\gamma\bar{\delta}}$, etc. From \eqref{2.28}, it follows that
\begin{align}
&f^\alpha_{k\bar{l}}=0,\label{2.30}\\
&f^\alpha_{kl}=f^\alpha_{lk}+\sum_if^\alpha_i\tau^i_{kl}-\sum_{\beta,\gamma}f^\beta_kf^\gamma_l\hat{\tau}^\alpha_{\beta\gamma},\\
&\sum_if^\alpha_i\tau^i_{\bar{k}\bar{l}}-\sum_{\beta,\gamma}f^{\bar{\beta}}_{\bar{k}}f^{\bar{\gamma}}_{\bar{l}}\hat{\tau}^\alpha_{\bar{\beta}\bar{\gamma}}=0.
\end{align}
Applying the exterior derivative to \eqref{2.29} and using \eqref{2.30} and the structure equations of $M$ and $\tilde{M}$ yield
\begin{align}
\sum_kDf^\alpha_{lk}\wedge\theta^k+\sum_jf^\alpha_j\Omega^j_l+\sum_kf^\alpha_{lk}\tau^k-\sum_\beta f^\beta_lf^*\Omega^\alpha_\beta=0,\label{2.33}
\end{align}
where
\begin{align}
Df^\alpha_{lk}=df^\alpha_{lk}-\sum_jf^\alpha_{lj}\theta^j_k-\sum_jf^\alpha_{jk}\theta^j_l+\sum_\beta f^\beta_{lk}\hat{\theta}^\alpha_\beta=\sum_m\left(f^\alpha_{lkm}\theta^m+f^\alpha_{lk\bar{m}}\theta^{\bar{m}}\right).\label{2.34}
\end{align}
Making use of \eqref{2.34}, we deduce from \eqref{2.33} that
\begin{align}
&f^\alpha_{lmn}=f^\alpha_{lnm}+\sum_jf^\alpha_jR^j_{lmn}+\sum_kf^\alpha_{lk}\tau^k_{mn}-\sum_{\beta,\gamma,\delta}\hat{R}^\alpha_{\beta\gamma\delta}f^\beta_lf^\gamma_mf^\delta_n,\\
&f^\alpha_{lm\bar{n}}=\sum_jf^\alpha_{j}R^j_{lm\bar{n}}-\sum_{\beta,\gamma,\delta}\hat{R}^\alpha_{\beta\gamma\bar{\delta}}f^\beta_lf^\gamma_mf^{\bar{\delta}}_{\bar{n}},\label{2.38.}\\
&\sum_{j}f^\alpha_jR^j_{l\bar{m}\bar{n}}+\sum_{k}f^\alpha_{lk}\tau^k_{\bar{m}\bar{n}}-\sum_{\beta,\gamma,\delta}\hat{R}^\alpha_{\beta\bar{\gamma}\bar{\delta}}f^\beta_lf^{\bar{\gamma}}_{\bar{m}}f^{\bar{\delta}}_{\bar{n}}=0.
\end{align}

Before ending this section, we give the definition of tamed exhaustion functions on almost Hermitian manifolds, which was introduced by Royden in \cite{[Ro]}.
\begin{definition}\label{def3.2}
Let $(M^m,J, g)$ be an almost Hermitian manifold. A continuous function $u: M\rightarrow \mathbb{R}$ is called a tamed exhaustion function, if it satisfies
\begin{enumerate}[(i)]
\item $u\geq 0$.
\item $u$ is proper, i.e., $u^{-1}((-\infty, c])$ is compact in $M$ for every constant $c\in \mathbb{R}$.
\item There exists a constant $C>0$ such that, for any $p\in M$, there is an open neighborhood $V_p$ of $p$ and a $C^2$ function $v: V_p\rightarrow \mathbb{R}$ with 
\begin{align}
v(x)\geq u(x)\ (\forall x\in V_p),\ v(p)=u(p),\ \|\nabla v\|(p)\leq C,\ (v_{k\bar{l}}(p))\leq C(g_{k\bar{l}}(p)).\label{2.49}
\end{align}
\end{enumerate}
Here the norm $\|\cdot\|$ is induced by the metric $g$, the matrix $(v_{k\bar{l}})$ is the Hessian of $v$ with respect to the canonical connection $\nabla$, $g_{k\bar{l}}$ are the components of almost Hermitian metric $g$, and the expression $A\leq B$ for two Hermitian matrices means that $B-A$ is positive semi-definite.
\end{definition}
We remark that the constant $C$ in Definition \ref{def3.2} is independent of the choice of the point $p\in M$.  However, for any point $p\in M$, there always exists a constant $C_p$ depending on $p$ and a $C^2$ function $v: V_p\rightarrow \mathbb{R}$ satisfying \eqref{2.49}.
\section{Tamed exhaustion functions on almost Hermitian manifolds}
In this section, we will establish the existence of tamed exhaustion functions on almost Hermitian manifolds. For this purpose, we need the following Hessian comparison theorem on almost Hermitian manifolds.

\begin{lemma}\label{lemma3.2}
Let $(M^m,J, g)$ be a complete almost Hermitian manifold with holomorphic bisectional curvature 
\begin{align}
R(X,\overline{X},Y,\overline{Y})\geq -B(1+r(x))^\alpha,\label{3.1}
\end{align}
the torsion 
\begin{align}
\|\tau(X,Y)\|\leq A_1(1+r(x))^\beta,\label{3.2}
\end{align}
and the $(2,0)$ part of the curvature tensor 
\begin{align}
|R(\overline{X},Y,Y,X)|\leq A_2(1+r(x))^\gamma,\label{3.3}
\end{align}
where $r(x)$ is the Riemannian distance of $M$ from a fixed point $x_0$ to $x\in M$, $X,Y\in T^{1,0}_{x}M$ with $\|X\|=\|Y\|=1$, and $B, A_1, A_2$ are positive constants, $\alpha\geq -2, \beta\geq 0,\gamma\geq 0$. Then
\begin{align}
\left(r_{i\bar{j}}\right)\leq \left\{\frac{1}{r}+\left[B(1+r)^\alpha+(4\sqrt{m}+3)A_1^2(1+r)^{2\beta}+2A_2(1+r)^\gamma \right]^{\frac{1}{2}}\right\}\left( g_{i\bar{j}}\right)
\end{align}
holds outside the cut locus of $x_0$.
\end{lemma}
\proof 
When $\alpha=\beta=\gamma=0$, this lemma has been proved in \cite{[Yu]}. For general cases, we will use a similar method to prove it. Let $\sigma(t)$ be a geodesic in $M$ with respect to $\nabla^{LC}$ with $\sigma(0)=x_0$ and $\|\sigma'(0)\|=1$. Let $\{e_i\}_{i=1}^m$ be a unitary frame field of $T^{1,0}M$ on a small neighborhood of $\sigma$. Since $\|\nabla r\|=1$ holds outside the cut locus of $x_0$, we have (cf. \cite[Lemma 3.2]{[Yu]})
\begin{align}
0&=\frac{1}{2}\|\nabla r\|^2_{k\bar{l}}\\
&=\left(\sum_ir_ir_{\bar{i}}\right)_{k\bar{l}}\notag\\
& =\sum_i\left(r_{k\bar{l}i}r_{\bar{i}}+r_{k\bar{l}\bar{i}}r_{i}\right)+\sum_ir_{k\bar{i}}r_{i\bar{l}}+\sum_{i,j}\left(r_{\bar{i}}\tau^j_{ik}r_{j\bar{l}}+r_{k\bar{j}}\tau^{\bar{j}}_{\bar{i}\bar{l}}r_i\right)\notag\\
&\ \ \ \ \ +\sum_ir_{ik}r_{\bar{i}\bar{l}}+\sum_{i,j}\left(r_{jk}\tau^j_{\bar{i}\bar{l}}r_i+r_{\bar{i}}\tau^{\bar{j}}_{ik}r_{\bar{j}\bar{l}}\right)+\sum_{i,j}
\left(R_{\bar{l}jik}+\sum_m\tau^{\bar{m}}_{ik}\tau^{\bar{j}}_{\bar{l}\bar{m}}\right)r_{\bar{j}}r_{\bar{i}}\notag\\
&\ \ \ \ \ \ \ \ +\sum_{i,j}\left(R_{k\bar{j}\bar{i}\bar{l}}+\sum_m\tau^j_{km}\tau^m_{\bar{i}\bar{l}}\right)r_jr_i+\sum_{i,j}\left(R_{i\bar{j}k\bar{l}}+\sum_m\left(\tau^{\bar{i}}_{km}\tau^{m}_{\bar{j}\bar{l}}+\tau^{\bar{m}}_{ik}\tau^{j}_{\bar{l}\bar{m}}\right)\right)r_jr_{\bar{i}}\notag.
\end{align}
Hence, 
\begin{align}
\frac{dX}{dr}+X^2+AX+XA^*&=-B^*B-B^*C-C^*B-(D+D^*)-E\\
&=-(B^*+C^*)(B+C)+C^*C-(D+D^*)-E\notag\\
&\leq C^*C-(D+D^*)-E,\notag
\end{align}
where matrices $X=(r_{k\bar{l}}), A=(\sum_jr_{\bar{j}}\tau^l_{jk}), B=(r_{\bar{k}\bar{l}}), C=(\sum_j\tau^k_{\bar{j}\bar{l}}r_j), D=(\sum_{i,j}(R_{k\bar{j}\bar{i}\bar{l}}+\sum_m\tau^j_{km}\tau^m_{\bar{i}\bar{l}})r_jr_i), E=(\sum_{i,j}(R_{i\bar{j}k\bar{l}}+\sum_m(\tau^{\bar{i}}_{km}\tau^{m}_{\bar{j}\bar{l}}+\tau^{\bar{m}}_{ik}\tau^{j}_{\bar{l}\bar{m}}))r_jr_{\bar{i}})$, and $B^*=\overline{B}^t, C^*=\overline{C}^t, D^*=\overline{D}^t$. Using the assumptions \eqref{3.1},\eqref{3.2},\eqref{3.3} and making a similar argument as in \cite[Theorem 2.1]{[Yu]} yield
\begin{align}
&\frac{dX}{dr}+X^2+AX+XA^*\\
&\ \ \ \ \ \ \ \  \leq \left(\frac{1}{2}B(1+r)^\alpha+\left(2\sqrt{m}+\frac{1}{2}\right)A_1^2(1+r)^{2\beta}+A_2(1+r)^\gamma\right)I_n,\notag\\
&X\leq \left(\frac{1}{r}+\frac{A_1}{\sqrt{2}}\right)I_n \ \text{as}\ r\rightarrow 0^+,\label{3.8}
\end{align}
where $I_n$ is the $n\times n$ identity matrix.
Set
\begin{align}
Y=\left(\frac{1}{r}+h(r)^{\frac{1}{2}} \right)I_n,
\end{align}
where 
\begin{align}
h(r)=B(1+r)^\alpha+(4\sqrt{m}+3)A_1^2(1+r)^{2\beta}+2A_2(1+r)^\gamma.\label{3.10}
\end{align}
Then by a simple computation, we obtain 
\begin{equation}
\frac{dY}{dr}\geq
\begin{cases}
-\frac{1}{r^2}I_n &\text{$\alpha\geq 0$}\\
\left(-\frac{1}{r^2}+\frac{\alpha B^{\frac{1}{2}}}{2(1+r)^{1-\frac{\alpha}{2}}}\right)I_n &\text{$-2\leq \alpha< 0$}
\end{cases}
\end{equation}
\begin{align}
&Y^2\geq \left(\frac{1}{r^2}+h(r)+\frac{(4\sqrt{m}+3)^{\frac{1}{2}}A_1(1+r)^{\beta}}{r}+\frac{B^{\frac{1}{2}}(1+r)^{\frac{\alpha}{2}}}{r}\right)I_n,\\
&AY+YA^*\geq -\sqrt{2}A_1(1+r)^\beta\left(\frac{1}{r}+h(r)^{\frac{1}{2}}\right)I_n\label{3.13}\\
&\ \ \ \ \ \ \ \ \ \ \ \ \ \ \geq \left(-\frac{\sqrt{2}A_1(1+r)^\beta}{r}-A_1^2(1+r)^{2\beta}-\frac{1}{2}h(r)\right)I_n,\notag
\end{align}
where \eqref{3.13} is due to 
\begin{align}
|u^*(A+A^*)u|&\leq |\sum_{j,k,l}u_{\bar{l}}\tau^l_{jk}r_{\bar{j}}u_k|+|\sum_{j,k,l}u_{\bar{l}}r_j\tau^{\bar{l}}_{\bar{j}\bar{k}}u_k|\\
&\leq \left(\sum_l|u_l|^2\right)^{\frac{1}{2}}\left\{\left(\sum_l|\sum_{j,k}\tau^l_{jk}r_{\bar{j}}u_k|^2\right)^{\frac{1}{2}}+\left(\sum_l |\sum_{j,k}\tau^{\bar{l}}_{\bar{j}\bar{k}}r_{\bar{j}}u_k|^2\right)^{\frac{1}{2}} \right\}\notag\\
&\leq \sqrt{2}A_1(1+r)^\beta\notag
\end{align}
for any unit column vector $u=(u_1,u_2,\dots, u_m)^t$, and 
\begin{align}
\sqrt{2}xy\leq x^2+\frac{1}{2}y^2.
\end{align}
Therefore, 
\begin{align}
\frac{dY}{dr}+&Y^2+AY+YA^*\\
&\geq \left(\frac{1}{2}B(1+r)^\alpha+\left(2\sqrt{m}+\frac{1}{2}\right)A_1^2(1+r)^{2\beta}+A_2(1+r)^\gamma\right)I_n\notag\\
&\geq \frac{dX}{dr}+X^2+AX+XA^*.\notag
\end{align}
From \eqref{3.8}-\eqref{3.10}, it follows that
\begin{align}
Y\geq X, \ \text{as\ $r\rightarrow 0^+$}.
\end{align}
In terms of the comparison theorem for matrix Ricatti equation (cf. \cite{[Ro-2]}), we get
\begin{align}
X\leq Y
\end{align}
outside the cut locus of $x_0$.
\qed
\begin{remark}
If $J$ is integrable, then the Nijenhuis tensor $\tau^i_{\bar{j}\bar{k}}\equiv0$ for any $i,j,k$. From \eqref{2.18.}, it follows that the $(2,0)$ part of the curvature tensor always vanishes, so the assumption \eqref{3.3} can be removed.
\end{remark}

Making use of the above lemma, we have the following 
\begin{theorem}\label{thm3.4}
Let $(M^m,J, g)$ be an almost Hermitian manifold as in Lemma \ref{lemma3.2} with $\alpha=\gamma=2, \beta=1$. Then it admits a tamed exhaustion function.
\end{theorem}
\proof
Consider the continuous function on $M$:
\begin{align}
u(x)=\log{(1+r(x)^2)},
\end{align}
where $r(x)=d(x_0,x)$ is the Riemannian distance from a fixed point $x_0\in M$ to a point $x$ in $M$. We will prove $u$ is a tamed exhaustion function on $M$. Since $M$ is complete, the function $u$ is proper. If $x\in M\setminus (Cut(x_0)\cup \{x_0\})$, where $Cut(x_0)$ denotes the cut locus of $x_0$, then the function $u$ is smooth near $x$ and satisfies 
\begin{align}
\|\nabla u\|(x)=\frac{2r(x)\|\nabla r\|}{1+r(x)^2}=\frac{2r(x)}{1+r(x)^2}\leq1.
\end{align}
Moreover, using Lemma \ref{lemma3.2} with $\alpha=\gamma=2, \beta=1$, we obtain
\begin{align}
\sum_{i,j}u_{i\bar{j}}(x)\xi^i\overline{\xi^j}&=\frac{2|\sum_ir_i\xi^i|^2}{1+r(x)^2}+\frac{2r(x)\sum_{i,j} r_{i\bar{j}}\xi^i\overline{\xi^j}}{1+r(x)^2}-\frac{4r(x)^2|\sum_k r_k\xi^k|^2}{(1+r(x)^2)^2}\label{3.21.}\\
&\leq\frac{2\left(\sum_i|r_i|^2\right)\left(\sum_i|\xi^i|^2\right)}{1+r(x)^2}+\frac{2r(x)\left(\frac{1}{r(x)}+h(r(x))^{\frac{1}{2}}\right)}{1+r(x)^2}\notag\\
&\leq\frac{3}{1+r(x)^2}+2C\frac{r(x)(1+r(x))}{1+r(x)^2}\notag
\end{align}
for any unite vector $\xi=(\xi^1,\xi^2,\dots,\xi^m)$, where $h(r)=B(1+r)^2+(4\sqrt{m}+3)A_1^2(1+r)^{2}+2A_2(1+r)^2$, and $C=[B+(4\sqrt{m}+3)A_1^2+2A_2]^{\frac{1}{2}}$. Since $\lim_{r\rightarrow 0+}\frac{(1+r)r}{1+r^2}=0$ and $\lim_{r\rightarrow +\infty}\frac{(1+r)r}{1+r^2}=1$, there exists a constant $C'$ independent of $x$ such that $\frac{(1+r)r}{1+r^2}<C'$, so $\sum_{i,j}u_{i\bar{j}}(x)\xi^i\overline{\xi^j}< 3+2CC'$.

If $x$ is a cut point of $x_0$, let $a$ be the distance between $x_0$ and the nearest cut-point of $x_0$, and take $x_1'\in B_{\frac{a}{2}}(x_0)\setminus\{x_0\}=\{x\in M: 0<r(x)<\frac{a}{2}\}$ on the minimal geodesic between $x_0$ and $x$, then $d(x_1',x)\geq \frac{a}{2}$. Set 
\begin{align}
v(z)=\log{\left(1+\left(r(x_1')+d(x_1',z)\right)^2\right)}.
\end{align}
for any $z$ in some small neighborhood $V_x$ of $x$. Clearly, ${v}\in C^\infty(V_x)$, ${v}(x)=u(x)$, ${v}(z)\geq u(z)$ for any $z\in V_x$. By a computation similar to \eqref{3.21.}, we obtain $\|\nabla {v}\|<1$ and 
\begin{align}
\sum_{i,j}{v}_{i\bar{j}}(x)\xi^i\overline{\xi^j}&\leq\frac{1}{1+r(x)^2}+\frac{2r(x)\left(\frac{1}{d(x_1',x)}+h(d(x_1',x))^{\frac{1}{2}}\right)}{1+r(x)^2}\\
&\leq1+\frac{1}{d(x_1',x)}+2C\frac{r(x)(1+r(x))}{1+r(x)^2}\notag\\
&\leq1+\frac{2}{a}+2CC'.\notag
\end{align}
Therefore, $u$ is a tamed exhaustion function of $M$.
\qed

Making a similar argument of \cite[Proposition 6.1]{[KL-2]} or \cite{[KL-1]}, it is easy to see that the Omori-Yau type maximum principle holds for any almost Hermitian manifold that admits a tamed exhaustion function. 
\begin{proposition}\label{prop3.5}
Let $(M^m,J, g)$ be an almost Hermitian manifold which admits a tamed exhaustion function. Then for any function $f:M\rightarrow \mathbb{R}$ with $\sup_M f<+\infty$, there exists a sequence $\{x_n\}\subset M$ such that
\begin{align}\label{3.26}
\lim_{n\rightarrow +\infty}f(x_n)=\sup_Mf,\ \lim_{n\rightarrow +\infty}\|\nabla f\|(x_n)=0,\ \limsup_{n\rightarrow +\infty} (\nabla df)(X_n,\overline{X_k})\leq 0
\end{align}
for any $X_n\in T^{1,0}_{x_n}M$ with $\|X_n\|=1$.
\end{proposition}

Combining Theorem \ref{thm3.4} and Proposition \ref{prop3.5}, we obtain the following corollary.
\begin{corollary}
Let $(M^m,J, g)$ be an almost Hermitian manifold as in Lemma \ref{lemma3.2} with $\alpha=\gamma=2, \beta=1$. Then for any function $f:M\rightarrow \mathbb{R}$ which is bounded above, there exists a sequence $\{x_n\}\subset M$ such that \eqref{3.26} holds.
\end{corollary}

\section{Schwarz type lemmas for almost Hermitian manifolds}
In this section, we will give the proofs of Theorem \ref{thm1.2} and Corollary \ref{cor1.6}.
\proof[Proof of Theorem \ref{thm1.2}]
At each $x\in M$, one can choose a unitary frame $\{\xi_i(x)\}$ for $T^{1,0}M$ such that 
\begin{align}
f^*\tilde{g}=\sum_{i,j,\alpha} f^\alpha_if^{\bar{\alpha}}_{\bar{j}}\theta^i\theta^{\bar{j}}=\sum_{i} \lambda_i(x) \theta^i\theta^{\bar{i}}
\end{align}
where $\lambda_1(x)\geq\lambda_2(x)\geq \cdots\geq \lambda_m(x)\geq 0$ are the eigenvalues of the positive semi-definite Hermitian matrix $\left(\sum_\alpha f^\alpha_if^{\bar{\alpha}}_{\bar{j}}\right)_x$. For simplicity, we write $\lambda_1(x)$ as $\lambda(x)$, then we have
\begin{align}
f^*\tilde{g}\leq \lambda(x) g,\label{4.2}
\end{align}
so it is sufficient to prove 
\begin{align}
\lambda(x)\leq\frac{k_1}{k_2}.
\end{align}
For this purpose, we consider the following function
\begin{align}
\phi(x)=(1-\epsilon u(x))^2\lambda(x)
\end{align}
for $x\in D_\epsilon=\{p\in M: u(p)<\frac{1}{\epsilon}\}$, where $u: M\rightarrow \mathbb{R}$ is a tamed exhaustion function. Since $u$ is nonnegative and proper, the closure $\overline{D_{\epsilon}}$ is compact. The nonnegative function $\phi$ is continuous on $\overline{D_{\epsilon}}$ and vanishes on the $\partial D_\epsilon$, hence $\phi$ attains its maximum at some point $x_\epsilon\in D_\epsilon$. In order to estimate the upper bound of $\lambda(x)$, we want to use the maximum principle for $\phi$ at $x_\epsilon$. However, both $u$ and $\lambda(x)$ are only continuous but may not be twice differentiable near $x_\epsilon$. This can be remedied by the following method: choose any local smooth unitary frame field $\{\eta_i\}$ for $T^{1,0}M$ such that $\eta_1(x_\epsilon)=\xi_1(x_\epsilon)$ is the eigenvector of $\lambda(x_\epsilon)$. Define
\begin{align}
\tilde{\lambda}=\|df(\eta_1)\|^2=\sum_\alpha f^\alpha_1f^{\bar{\alpha}}_{\bar{1}}.
\end{align}
 It is obvious that $\tilde{\lambda}$ is smooth near $x_\epsilon$ and 
 \begin{align}
 \tilde{\lambda}(x)\leq \lambda(x),\ \ \ \tilde{\lambda}(x_\epsilon)=\lambda(x_\epsilon).
 \end{align}
 Let 
 \begin{align}
 \tilde{\phi}(x)=(1-\epsilon v(x))^2\tilde{\lambda}(x),
 \end{align}
 where $v: V_{x_\epsilon}\rightarrow \mathbb{R}$ is the $C^2$ function defined as in Definition \ref{def3.2}. Then we obtain
 \begin{align}
 \tilde{\phi}(x)\leq (1-\epsilon u(x))^2\lambda(x)\leq  (1-\epsilon u(x_\epsilon))^2\lambda(x_\epsilon)=\tilde{\phi}(x_\epsilon)
 \end{align}
 for any $x\in V_{x_\epsilon}$, hence $\tilde{\phi}(x)$ has a local maximum point $x_\epsilon$. Applying the maximum principle for $\tilde{\phi}(x)$ at $x_\epsilon$, we deduce that
 \begin{equation}
 \nabla \tilde{\lambda}=\frac{2\epsilon\tilde{\lambda}}{1-\epsilon v}\nabla v,\label{411}
 \end{equation}
 and 
 \begin{align}
 0\geq &(\tilde{\phi})_{1\bar{1}}\label{412}\\
 =&(1-\epsilon v)^2{\tilde{\lambda}}_{1\bar{1}}-2\epsilon(1-\epsilon v)(\tilde{\lambda}_1v_{\bar{1}}+\tilde{\lambda}_{\bar{1}}v_1)\notag\\
 &+2\epsilon^2\tilde{\lambda}|v_1|^2-2\epsilon(1-\epsilon v)\tilde{\lambda}v_{1\bar{1}}.\notag
 \end{align}
 Substituting \eqref{411} into \eqref{412} yields
 \begin{align}
 0\geq (1-\epsilon v)^2{\tilde{\lambda}}_{1\bar{1}}-6\epsilon^2\tilde{\lambda}|v_1|^2-2\epsilon\tilde{\lambda}(1-\epsilon v)v_{1\bar{1}}.
 \end{align}
 From Definition \ref{def3.2}, we can see that
 \begin{align}
 0\leq |v_1|^2\leq \frac{1}{2}|\nabla v|^2\leq \frac{1}{2}C^2,\ \ \ v_{1\bar{1}}\leq Cg_{1\bar{1}}=C,
 \end{align}
 where $C$ is a constant defined as in Definition \ref{def3.2}. Hence,
\begin{align}
0\geq (1-\epsilon v)^2{\tilde{\lambda}}_{1\bar{1}}-3C^2\epsilon^2\tilde{\lambda}-2C\epsilon\tilde{\lambda}.\label{4.13}
\end{align}
By \eqref{2.30} and \eqref{2.38.}, we perform the following computation:
\begin{align}
\tilde{\lambda}_{1\bar{1}}&=\left(\sum_\alpha f^{\alpha}_1f^{\bar{\alpha}}_{\bar{1}}\right)_{1\bar{1}}\label{4.15}\\
&=\left(\sum_\alpha f^{\alpha}_{11}f^{\bar{\alpha}}_{\bar{1}}\right)_{\bar{1}}\notag\\
&=\sum_\alpha f^{\alpha}_{11\bar{1}}f^{\bar{\alpha}}_{\bar{1}}+\sum_\alpha |f^\alpha_{11}|^2\notag\\
&\geq-\sum_{\alpha,\beta,\gamma,\delta}f^{\bar{\alpha}}_{\bar{1}}f^\beta_1f^\gamma_1f^{\bar{\delta}}_{\bar{1}}\hat{R}^\alpha_{\beta\gamma\bar{\delta}}+\sum_jR^j_{11\bar{1}}\left(\sum_\alpha f^\alpha_jf^{\bar{\alpha}}_{\bar{1}}\right)\notag
\end{align}
In terms of the curvature assumptions for both $M$ and $N$ and $\left(\sum_\alpha f^\alpha_jf^{\bar{\alpha}}_{\bar{1}}\right)(x_\epsilon)=\delta_{1j}\tilde{\lambda}(x_\epsilon)$, we have, at $x_\epsilon$,
\begin{align}
\tilde{\lambda}_{1\bar{1}}\geq k_2\tilde{\lambda}^2-k_1\tilde{\lambda}.\label{4.16}
\end{align}
From \eqref{4.13} and \eqref{4.16}, we deduce that
\begin{align}
(1-\epsilon u(x))^2\lambda(x)&\leq (1-\epsilon u(x_\epsilon))^2\lambda(x_\epsilon)=(1-\epsilon v(x_\epsilon))^2\tilde{\lambda}(x_\epsilon)\\
&\leq \frac{k_1}{k_2}+\frac{3C^2\epsilon^2+2C\epsilon}{k_2}\notag
\end{align}
for any $x\in D_\epsilon$. Let $\epsilon\rightarrow 0+$, we obtain
\begin{align}
\lambda(x)\leq \frac{k_1}{k_2}
\end{align}
for any $x\in M$. 
\qed

When $M$ is compact, similar to the proof of Theorem \ref{thm1.2}, we can give a simple proof of Corollary \ref{cor1.5}.
\proof[Proof of Corollary \ref{cor1.5}] Consider the largest eigenvalue $\lambda(x)$ defined as in \eqref{4.2} which is continuous on $M$. Assume that $f$ is not constant, then $\lambda(x)\not\equiv 0$. Since $M$ is compact, there exists a point $x_0\in M$ such that $\lambda(x_0)=\max_M\lambda>0$. Choose any local smooth unitary frame field $\{\eta_i\}$ for $T^{1,0}M$ such that $\eta_1(x_0)$ is the eigenvector of $\lambda(x_0)$. Define
\begin{align}
\tilde{\lambda}=\|df(\eta_1)\|^2=\sum_\alpha f^\alpha_1f^{\bar{\alpha}}_{\bar{1}}.
\end{align}
 It is obvious that $\tilde{\lambda}$ is smooth near $x_0$ and 
 \begin{align}
 \tilde{\lambda}(x)\leq \lambda(x),\ \ \ \tilde{\lambda}(x_0)=\lambda(x_0).
 \end{align}
Therefore, $x_0$ is also the local maximum point of $\tilde{\lambda}$. According to the Hessian type Bochner formula \eqref{4.15} and the maximum principle, we have, at $x_0$
\begin{align}
0\geq \tilde{\lambda}_{1\bar{1}}&\geq-\sum_{\alpha,\beta,\gamma, \delta}f^{\bar{\alpha}}_{\bar{1}}f^\beta_1f^\gamma_1f^{\bar{\delta}}_{\bar{1}}\hat{R}^\alpha_{\beta\gamma\bar{\delta}}+R^1_{11\bar{1}}\tilde{\lambda}>0,
\end{align}
 which leads to a contradiction. Therefore, $\lambda\equiv 0$, i.e., $f$ is constant.
 \qed

\bigskip
\bigskip

Weike Yu

School of Science 

Nanjing University of Science and Technology

Nanjing, 210094, Jiangsu, P. R. China

wkyu2018@outlook.com

\bigskip

\end{document}